\documentclass[12pt]{amsart}
\usepackage{amssymb}
\usepackage{amscd}
\newcommand{\nc}{\newcommand}

\nc{\tB}{\widetilde{\CB}}
\nc{\hH}{\widehat{H}}
\nc{\BI}{\CB_{\infty}}
\nc{\CGI}{\CG_{\infty}}
\nc{\LBA}{\mathtt{LBA}_0}
\nc{\Locc}{\mathtt{Locc}}
\nc{\Liealg}{\mathtt{Lie}}
\nc{\HA}{\mathtt{HA}_0}

\nc{\Del}{\operatorname{Del}}
\nc{\Ind}{\operatorname{Ind}}
\nc{\OT}{\mathbb{T}}
\nc{\ad}{\operatorname{ad}}
\nc{\Vect}{\mathtt{Vect}}
\nc{\Vectgr}{\mathtt{Vectgr}}
\nc{\Endop}{\mathrm{Endop}}
\nc{\Har}{\mathrm{Harr}}

\nc{\invtensor}{\underset{\leftarrow}{\otimes}}
\nc{\rlarrows}{\begin{picture}(1,0.4)
                \put(0,-0.1){\makebox(1,0.2){$\leftarrow$}}
                \put(0,0.1){\makebox(1,0.2){$\to$}}
              \end{picture}}
\nc{\rra}{\begin{picture}(1,0.4)
                \put(0,-0.1){\makebox(1,0.2){$\lra$}}
                \put(0,0.1){\makebox(1,0.2){$\lra$}}
              \end{picture}}
\nc{\Left}{\mathbf L}  
\nc{\Right}{\mathbf R} 
\nc{\gr}{\operatorname{gr}}
\nc{\Ho}{\operatorname{Ho}}
\nc{\alt}{\operatorname{alt}}
\nc{\Sym}{\operatorname{Sym}}
\nc{\sym}{\operatorname{sym}}
\nc{\id}{\operatorname{id}}
\nc{\Der}{\operatorname{Der}}
\nc{\im}{\operatorname{Im}}
\nc{\Ker}{\operatorname{Ker}}
\nc{\coker}{\operatorname{Coker}}
\nc{\Col}{\operatorname{Col}}
\nc{\ter}{\operatorname{ter}}
\nc{\intl}{\operatorname{int}}
\nc{\val}{\operatorname{val}}

\nc{\TN}{{\cal N}}
\nc{\Norm}{\operatorname{N}}
\nc{\Nor}{\operatorname{N}}
\nc{\Tor}{\operatorname{Tor}}
\nc{\res}{\operatorname{res}}
\nc{\Stab}{\operatorname{Stab}}
\nc{\Hom}{\operatorname{Hom}}
\nc{\chom}{\CH\!o\!m}
\nc{\uhom}{\CH\!o\!m}
\nc{\End}{\operatorname{End}}
\nc{\holim}{\operatorname{holim}}
\nc{\dirlim}{\underset{\rightarrow}{\lim}\,}
\nc{\invlim}{\underset{\leftarrow}{\lim}\,}
\nc{\com}{\operatorname{co}}
\nc{\Tot}{\operatorname{Tot}}
\nc{\Th}{\operatorname{Th}}
\nc{\Cech}{\check{C}}
\nc{\Spec}{\operatorname{Spec}}
\nc{\Spf}{\operatorname{Spf}}
\nc{\MC}{\operatorname{MC}}
\nc{\U}{\operatorname{U}}
\nc{\Diff}{{\cal D}\mbox{\em iff}}
\nc{\Mor} {{\cal M}or}
\nc{\Ob}{\operatorname{Ob}}
\nc{\cone}{\widehat}
\nc{\Coder}{\operatorname{Coder}}
\nc{\pr}{\operatorname{pr}}
\nc{\diag}{\operatorname{diag}}

\nc{\CHo}{{\cal{H}\mbox{\it{o}}}}
\nc{\Mod}{{\mathtt{mod}}}       
\nc{\Modf}{{\mathtt{modf}}}       
\nc{\Modg}{{\mathtt{modg}}}       
\nc{\Ab}{{\mathtt {Ab}}}          
\nc{\Alg}{{\mathtt {Alg}}} 
\nc{\Hoalg}{{\mathtt {Hoalg}}} 
\nc{\Valg}{{\mathtt {Viral}}} 
\nc{\Algf}{{\mathtt {Algf}}} 
\nc{\Algg}{{\mathtt {Algg}}} 
\nc{\Coalg}{{\mathtt {Coalg}}} 
         
\nc{\dgc}{{\mathtt{dgc}}}
\nc{\dgca}{{\mathtt{dgca}}}
\nc{\dgcu}{{\mathtt{dgcu}}}
\nc{\dgcuf}{{\mathtt{dgcuf}}}
\nc{\dgcf}{{\mathtt{dgcf}}}
\nc{\dgcg}{{\mathtt{dgcg}}}
\nc{\dgcc}{{\mathtt{dgccc}}}

\nc{\dgl}{{\mathtt{dglie}}}
\nc{\dgla}{{\mathtt{dgla}}}
\nc{\dglf}{{\mathtt{dglf}}}
\nc{\dglg}{{\mathtt{dglg}}}
\nc{\dga}{{\mathtt{dga}}}
\nc{\art}{{\mathtt {art}}}
\nc{\dgar}{{\mathtt {dgart}^{\leq 0}}}
\nc{\simpl}{\Delta^{\op}\Ens}
\nc{\Coll}{{\mathtt{Coll}}}
\nc{\Kan}{{\mathtt {Kan}}}

\nc{\Grp}{{\mathtt {Grp}}}
\nc{\Cat}{{\mathtt {Cat}}}
\nc{\Ens}{{\mathtt {Ens}}}
\nc{\op}{{\operatorname{op}}}
\nc{\Op}{{\mathtt{Op}}}

\nc{\Lie}{{\mathtt{LIE}}}
\nc{\Com}{{\mathtt{COM}}}
\nc{\Ass}{{\mathtt{ASS}}}

\nc{\pa}{\partial}

\nc{\cal}{\mathcal} 

\nc{\CA}{\cal A}
\nc{\CBB}{\cal B}
\nc{\CB}{\cal B}
\nc{\CC}{\cal C}
\nc{\CDD}{\cal D}
\nc{\CE}{\cal E}
\nc{\CF}{\cal F}
\nc{\CG}{\cal G}
\nc{\CH}{\cal H}
\nc{\CI}{\cal I}
\nc{\CJ}{\cal J}
\nc{\CK}{\cal K}
\nc{\CL}{\cal L}
\nc{\CM}{\cal M}
\nc{\CN}{\cal N}
\nc{\CO}{\cal O}
\nc{\CP}{\cal P}
\nc{\CQ}{\cal Q}
\nc{\CR}{\cal R}
\nc{\CS}{\cal S}
\nc{\CT}{\cal T}
\nc{\CU}{\cal U}
\nc{\CV}{\cal V}
\nc{\CW}{\cal W}
\nc{\CZ}{\cal Z}

\nc{\fa}{\mathfrak a}
\nc{\fg}{\mathfrak g}
\nc{\fk}{\mathfrak k}
\nc{\fh}{\mathfrak h}
\nc{\fm}{\mathfrak m}
\nc{\fn}{\mathfrak n}
\nc{\fA}{\mathfrak A}
\nc{\fC}{\mathfrak C}
\nc{\fI}{\mathfrak I}
\nc{\fS}{\mathfrak S}

\nc{\nen}{\newenvironment}
\nc{\ol}{\overline}
\nc{\ul}{\underline}
\nc{\lra}{\longrightarrow}
\nc{\lla}{\longleftarrow}
\nc{\Lra}{\Longrightarrow}
\nc{\Lla}{\Longleftarrow}
\nc{\Llra}{\Longleftrightarrow}
\nc{\hra}{\hookrightarrow}
\nc{\iso}{\overset{\sim}{\lra}}

\nc{\notebox}[1]{\noindent\fbox{\parbox{12.5cm}{\sf #1}}\\[8pt]}
\nc{\Thm}[1]{Theorem~\ref{#1}}
\nc{\Prop}[1]{Proposition~\ref{#1}}
\nc{\Lem}[1]{Lemma~\ref{#1}}
\nc{\Cor}[1]{Corollary~\ref{#1}}
\nc{\Conj}[1]{Conjecture~\ref{#1}}
\nc{\Claim}[1]{Claim~\ref{#1}}
\nc{\Defn}[1]{Definition~\ref{#1}}
\nc{\Exa}[1]{Example~\ref{#1}}
\nc{\Rem}[1]{Remark~\ref{#1}}
\nc{\Note}[1]{Note~\ref{#1}}

\nen{thm}[1]{\label{#1}{\bf Theorem.\ } \em}{}
\nen{prop}[1]{\label{#1}{\bf Proposition.\ } \em}{}
\nen{lem}[1]{\label{#1}{\bf Lemma.\ } \em}{}
\nen{klem}[1]{\label{#1}{\bf Key Lemma.\ } \em}{}
\nen{cor}[1]{\label{#1}{\bf Corollary.\ } \em}{}
\nen{conj}[1]{\label{#1}{\bf Conjecture.\ } \em}{}
\nen{claim}[1]{\label{#1}{\bf Claim.\ } \em}{}

\nen{defn}[1]{\label{#1}{\bf Definition.\ } }{}
\nen{exa}[1]{\label{#1}{\bf Example.\ } }{}

\nen{rem}[1]{\label{#1}{\em Remark.\ } }{}
\nen{note}[1]{\label{#1}{\em Note.\ } }{}
\nen{exer}[1]{\label{#1}{\em Exercise.\ } }{}
\nen{pf}{\begin{proof}}{\end{proof}}

\begin{document}

\title[Tamarkin's proof]
{Tamarkin's proof of Kontsevich formality theorem}
\author{Vladimir Hinich}
\address{Department of Mathematics, University of Haifa,
Mount Carmel, Haifa 31905,  Israel}
\email{hinich@math.haifa.ac.il}
\maketitle

\section{Introduction}

\subsection{}
This is an extended version of lectures given at Luminy colloquium
``Quantification par d\'eformation'' held  at December, 1999.

In this note we explain Tamarkin's proof \cite{t} of 
the following affine algebraic version of Kontsevich's formality theorem.

\subsection{}
\begin{thm}{KFT}Let $A$ be a polynomial algebra over a field $k$ of 
characteristic zero and let $\CC=C^*(A;A)$ be the cohomological Hochschild
complex of $A$ with coefficients at $A$. The dg Lie algebra $\CC[1]$
is formal, that is $\CC[1]$ is isomorphic in the homotopy category
of dg Lie algebras to its cohomology.
\end{thm}

Our sources are the original Tamarkin's paper~\cite{t} and the recent
preprint of Tamarkin-Tsygan~\cite{tt} where a simplification
of the original proof is sketched.

We tried to provide all necessary details that were sometimes
difficult to find in~\cite{t}. 

\subsection{}
In the first part (Sections~\ref{basic}--~\ref{kd}) we review some basic 
facts on operads and Koszul operads. In Section ~\ref{df} we study formality
of algebras over a Koszul operad. Following Halperin-Stasheff~\cite{hs},
we call a graded algebra $H$ over a Koszul operad $\CO$ {\em intrinsically
formal} if any dg $\CO$-algebra with cohomology isomorphic to $H$ is formal. 
We prove \Thm{crit-i-f} which gives a sufficient condition of intrinsic 
formality of a graded algebra over a Koszul operad in terms of its cohomology.

In Section~\ref{hoch} we calculate the cohomology
of Gerstenhaber algebra $H(\CC)$, $\CC$ being the Hochschild complex of a 
smooth $k$-algebra. The calculation shows that $H(\CC)$ is intrinsically
formal for a polynomial algebra. 
This proves that Kontsevich's formality theorem follows from
a version of Deligne's conjecture~\ref{dc} on the existence of 
homotopy Gerstenhaber algebra structure on the Hochschild complex. 
\Thm{dc} is proven in Sections~\ref{bet} and~\ref{two}.

{\em Acknowledgement.} I am grateful to B.~Tsygan for sending early 
version of his manuscript~\cite{tt} and to D.~Tamarkin for an
important remark. This note was written during
my stay at MSRI. I am grateful to MSRI for excellent working conditions.

\section{Basic definitions}
\label{basic}

In this section we recall the basic definitions of operads and operad
algebras. The most convenient reference here is ~\cite{gj}, ch.~1.

\subsection{Operads}
Let $\Vect$ be the category of vector spaces over a field $k$ of 
characteristic zero.

By definition, an $\mathbb{S}$-object in $\Vect$ is a collection 
$X=\{X(n)\},\ n\ge 0,$ of objects of $\Vect$ endowed with a right
action of the symmetric groups $S_n$. The category of $\mathbb{S}$-objects in
$\Vect$ admits a (non-symmetric) monoidal structure defined as follows.

Any $\mathbb{S}$-object $X$ defines a functor $\CS(X)$ (Schur functor) on 
$\Vect$ by the formula
\begin{equation}
\label{S(X)}
\CS(X): V\mapsto\bigoplus X(n)\otimes_{S_n}V^{\otimes n}.
\end{equation}

The monoidal operation on the category of $\mathbb{S}$-vector spaces 
is uniquely defined by the property
$$\CS(X\circ Y)=\CS(X)\circ\CS(Y).$$

\subsubsection{}
\begin{defn}{operad}An operad $\CO=\{\CO(n)\}$ in $\Vect$ is a monoid in the 
category of $\mathbb{S}$-vector spaces. The category of operads in $\Vect$ is
denoted $\Op(\Vect)$.
\end{defn}

In more conventional terms, an operad is an $\mathbb{S}$-vector space
$\{\CO(n)\}$ endowed with equivariant operations
\begin{equation}
\label{op-pro}
\CO(n)\otimes\CO(m_1)\otimes\ldots\otimes\CO(m_n)\to\CO(\sum m_i)
\end{equation}
and with a unit element $1\in\CO(1)$ satisfying natural associativity
and unit conditions.

\subsubsection{}For any vector space $V$ one defines an operad
$\Endop(V)$ to be a $\mathbb{S}$-vector space 
$$n\mapsto\Hom(V^{\otimes n}, V)$$
with the obvious composition and action of the symmetric groups.

\subsubsection{}
\begin{defn}{alg}
An algebra $A$ over an operad $\CO$ is a map of operads
$$ \CO\to\Endop(A).$$
\end{defn}
In other terms, an $\CO$-algebra structure on $A$ is given by a collection
 of $S_n$-equivariant maps 
$$ \CO(n)\otimes A^{\otimes n}\to A$$
satisfying natural associativity and unit properties.

\subsubsection{Examples}There are operads $\Ass$, $\Com$, $\Lie$ such that
corresponding algebras are associative, commutative and Lie algebras
respectively.

\subsection{Other tensor categories}

The definitions of the previous subsection make sense in any tensor 
(= monoidal symmetric) category $\CA$. The following cases will be of a special
interest for us.

\subsubsection{} $\CA=\Vectgr$ --- the category of $\mathbb{Z}$-graded
vector spaces. The commutativity constraint $X\otimes Y\iso Y\otimes X$
is defined by the standard formula
\begin{equation}
\label{signs}
x\otimes y\mapsto (-1)^{|x||y|}y\otimes x.
\end{equation}

\subsubsection{} $\CA=C(k)$ --- the category of complexes over $k$. The 
commutativity constraint in this case is given by the same 
formula~(\ref{signs}).

\subsubsection{}Let $\CO\in\Op(\CA)$ for a tensor category $\CA$ and let
$\alpha:\CA\to\CB$ be a tensor functor. Then $\alpha(\CO)$ is an operad
over $\CB$. This obvious construction allows one, for example, to consider
graded or dg Lie algebras as algebras over the operad $\Lie$ in $\Vectgr$ or
in $C(k)$ respectively.
 
\subsubsection{}Denote $k[n]$ to be a standard one-dimensional vector
space concentrated at degree $-n$. One defines the $n$-shift functor
$X\mapsto X[n]$ by the formula
$$ X[n]=k[n]\otimes X.$$
This formula makes sense both in $\Vectgr$ and in $C(k)$.

Let $\CO$ be an operad in $\Vectgr$ or $C(k)$. There is a uniquely
defined operad $\CO\{m\}$ such that a $\CO\{m\}$-algebra structure on $X$
is equivalent to a $\CO$-algebra structure on $X[m]$. One has
$$ \CO\{m\}(n)=\Lambda_n^{\otimes m}\otimes\CO(n)$$
where $\Lambda_n$ denotes the graded vector space (or complex) $k[n-1]$ 
endowed with the sign representation of the symmetric group $S_n$.

\subsection{Free algebras and free operads}

\subsubsection{}
Let $\CO$ be an operad over a tensor category $\CA$. Let $V$
be an $\mathbb{S}$-object in $\CA$.
The free $\CO$-algebra generated by $V$ is defined to be
\begin{equation}
\label{free}
\mathbb{F}_{\CO}(V)=\bigoplus_{n\ge 0}\CO(n)\otimes_{S_n}V^{\otimes n}
\end{equation}
with a canonical $\CO$-algebra structure. 

\subsubsection{} Let $X$ be an $\mathbb{S}$-object in $\CA$. The forgetful
functor from the category of operads to the category of $\mathbb{S}$-objects
in $\CA$ admits a left adjoint {\em free operad functor.}
Free operad $\mathbb{T}(X)$ generated by $X$ has an explicit description 
as a direct sum over trees (see~\cite{gj}, 1.4). 

\subsection{Cooperads and coalgebras}
\subsubsection{}
The notions of operad and algebra can be dualized. Let $\CC$ be a cooperad.
A $\CC$-coalgebra $X$ is called {\em nilpotent} if 
\begin{equation}
\label{nilpotent}
X=\cup_n\Ker(X\to\CC(n)\otimes X^{\otimes n}).
\end{equation}

From now on all coalgebras will be supposed to be nilpotent. 
We define $\Coalg(\CC)$ to be the category of nilpotent
$\CC$-coalgebras.  If $V$ is an $\mathbb{S}$-object in $\CA$, 
the cofree (nilpotent) cooalgebra cogenerated by $V$ is defined to be
\begin{equation}
\label{cofree}
\mathbb{F}^*_{\CC}(V)=\bigoplus_{n\ge 0}\left(\CC(n)\otimes V^{\otimes n}
\right)^{S_n}.
\end{equation}

Let $X$ be a $\CC$-coalgebra and $V$ be an $\mathbb{S}$-object. Any
map $X\to V$ of $\mathbb{S}$-objects defines canonically a map
of $\CC$-coalgebras $X\to\mathbb{F}^*_{\CC}(V)$. $V$ is called
an $\mathbb{S}$-object of {\em cogenerators} if the above map is injective.

Cofree cooperad cogenerated by $V$ is denoted $\mathbb{T}^*(V)$.
It is isomorphic to $\mathbb{T}(V)$ as an $\mathbb{S}$-object. However, we
prefer to have a different notation to stress that this is a cooperad.

\subsubsection{}
Let $\CO\in\Op(\Vectgr)$ be an operad such that $\CO(n)$ are all finite 
dimensional. Then the collection $\{\CO(n)^*\}$ admits an obvious structure 
of cooperad. This cooperad is denoted by $\CO^*$. 
Coalgebras over $\CO^*$ are sometimes called $\CO$-coalgebras. In the same
style, we will sometimes write $\mathbb{F}^*_{\CO}(V)$ instead of
$\mathbb{F}^*_{\CO^*}(V)$. Thus,  $\Com$-coalgebras are just cocommutative
coalgebras, $\Lie$-coalgebras are Lie coalgebras, etc.

\section{Koszul duality}
\label{kd}

\subsection{Quadratic operads and quadratic duals}
\subsubsection{}
\begin{defn}{}
An operad $\CO$ of graded vector spaces is called {\em quadratic} if it is
generated (as operad) by $\CO(2)$ and has only relations of valence $3$.

The latter condition means the following. 
Let  $V$ be the $\mathbb{S}$-object in $\Vectgr$ defined by the properties
$V(2)=\CO(2),\quad V(n)=0$ for $n\ne 2$. Since $\CO$ is generated by
its binary operations, the natural map $\mathbb{T}(V)\to\CO$ is surjective.

The operad $\CO$ is quadratic if the kernel of this map is generated 
(as an ideal in an operad) by an $S_3$-invariant subspace 
$R\subseteq\mathbb{T}(V)(3)$.
\end{defn}

Note that $\mathbb{T}(V)(3)=\Ind_{S_2}^{S_3}(V\otimes V)$ where $S_2$ acts on
the tensor product through the trivial action on the first factor.

A quadratic operad $\CO$ with generators $V$ and relations $R$ can be
described as the pushout in the category of operads
$$
\begin{CD}
\mathbb{T}(R)@>>>\mathbb{T}(V)\\
@VVV              @VVV        \\
* @>>>\CO
\end{CD}
$$
where $*$ denotes the trivial operad 
$$*(1)=k,\quad *(n)=0\text{ for }n\ne 1.$$

Dually, let $V$ be a graded vector space endowed with an action of $S_2$
and let $R$ be an $S_3$-invariant subspace of $\mathbb{T}^*(V)(3)$.
Denote $Q=\mathbb{T}^*(V)(3)/R$. Then a quadratic cooperad $\CC$ cogenerated
by $V$ with co-relations $R$ is defined as the pullback
$$
\begin{CD}
\CC @>>>\mathbb{T}^*(V)\\
@VVV              @VVV        \\
* @>>>\mathbb{T}^*(Q)
\end{CD}
$$

\subsubsection{}
\begin{defn}{}
\fbox{rewrite!!!}
1. Let $\CO$ be a quadratic operad with $V=\CO(2)$ and the space of relations
$R$. The dual cooperad $\CO^{\perp}$ is cogenerated by the space $V[1]$
with co-relations $\CO^{\perp}(3)=R[2]$.

2. Dually, for a cooperad $\CC$ cogenerated by $V$ with co-relations $Q$,
the quadratic dual operad $\CC^{\perp}$ is generated by $V[-1]$ with relations
given by the kernel
$$ \Ker(Q[-2]\to (V[-1]\circ V[-1])(3).$$
\end{defn}

\subsubsection{}
{\em Examples}

The operads $\Com$, $\Ass$, $\Lie$ are quadratic. Their quadratic dual
cooperads are given by the formulas
\begin{itemize}
\item $\Com^{\perp}=(\Lie\{-1\})^*,$ 
\item $\Ass^{\perp}=(\Ass\{-1\})^*,$
\item $\Lie^{\perp}=(\Com\{-1\})^*.$
\end{itemize}

\subsubsection{}
\begin{defn}{o-infty}
Let $\CO$ be a (graded) quadratic operad. A structure of $\CO_{\infty}$-algebra
on $X\in\Vectgr$ is given by a differential on the 
cofree $\CO^{\perp}$-coalgebra cogenerated by $X$. 
\end{defn}

\subsubsection{}
The above definition gives rise to an operad $\CO_{\infty}$
{\em in the category of complexes $C(k)$}.

Let $X$ have a structure of $\CO_{\infty}$-algebra. The differential
\begin{equation}
\label{diff}
Q:\mathbb{F}^*_{\CO^{\perp}}(X)\to\mathbb{F}^*_{\CO^{\perp}}(X)[1]
\end{equation}
is defined uniquely by its composition with the projection onto
the degree one component $F^{*1}_{\CO^{\perp}}(X)=X$. Thus, the
differential is given by the collection of maps
\begin{equation}
\label{diff-i}
Q_i:\mathbb{F}^{*i}_{\CO^{\perp}}(X)=\left(\CO^{\perp}(i)
\otimes X^{\otimes i}\right)^{S_i}\to X[1].
\end{equation}
in particular, $d:=Q_1$ defines a differential on $X\in\Vectgr$.

Define $\CO_{\infty}=\mathbb{T}(\CO^{\perp})$ to be the free
graded operad generated by $\CO^{\perp}$. The collection of maps
$Q_i$ from~(\ref{diff-i}) defines an action of $\CO_{\infty}$ on $X$. One
endows $\CO_{\infty}$ with a differential so that the condition $Q^2=0$
is equivalent to the statement that the action of $\CO_{\infty}$ on 
$(X,d=Q_1)$ respects the differentials.

\subsubsection{}
\begin{exa}{strict}
Let $X$ be a complex endowed with a $\CO$-algebra structure (dg $\CO$-algebra).
Define the differential $Q$ on $\mathbb{F}^*_{\CO^{\perp}}(X)$ as follows.

$Q_1:X\to X[1]$ is the differential of $X$. $Q_2:\CO^{\perp}(2)\otimes 
X^{\otimes 2}\to X[1]$ is defined by the $\CO$-algebra structure on $X$
since $\CO^{\perp}(2)=\CO(2)[1]$. $Q_i$ are defined to be zero for $i>2$.

The condition $Q^2=0$ can be easily verified. This means that any 
$\CO$-algebra admits a canonical $\CO_{\infty}$-algebra structure. 
\end{exa}

\Exa{strict} shows there is a canonical map of operads in $C(k)$
\begin{equation}
\label{strictishom}
\CO_{\infty}\to\CO
\end{equation}
(Here $\CO$ is supposed to have zero differential).

\subsubsection{}
\begin{defn}{koszul}
A quadratic operad  $\CO$ is called Koszul if the natural 
map~(\ref{strictishom}) is a quasi-isomorphism.
\end{defn}

Let $\CO$ be a quadratic operad and let $X$ be an $\CO_{\infty}$-algebra
(for instance, an $\CO$-algebra). The homology of $X$, $H_{\CO}(X)$,
is defined to be the homology of the complex 
$(\mathbb{F}^*_{\CO^{\perp}}(X),Q)$.

If $X=\mathbb{F}_{\CO}(V)$ for a graded vector space $V$, one has a canonical
map of complexes
\begin{equation}
\label{for-free}
(\mathbb{F}^*_{\CO^{\perp}}(X),Q)\to V.
\end{equation}

The following result can be used to prove koszulity of a quadratic operad.
\subsubsection{}
\begin{thm}{koszul-crit}(cf.~\cite{gk}, Thm.~4.2.5) A quadratic operad $\CO$
is Koszul iff for any graded vector space $V$ the canonical 
map~(\ref{for-free}) is quasi-isomorphism.
\end{thm}

\Thm{koszul-crit} implies that the operads $\Com,\Ass,\Lie$ are Koszul.

\section{Deformations and formality}
\label{df}

\subsection{Intrinsic formality}
In this section $\CO$ is a fixed Koszul operad.

\subsubsection{} 
\begin{defn}{formal}
A $\CO_{\infty}$-algebra $X$
is called to be {\em formal} if there exists a pair of quasi-isomorphisms
of $\CO_{\infty}$-algebras $X\leftarrow F\to H(X)$.
\end{defn}

\subsubsection{} 
\begin{defn}{intr-formal}
A graded $\CO$-algebra $H$ is {\em intrinsically formal} if any
$\CO_{\infty}$-algebra $X$ with $H(X)=H$ is formal.
\end{defn}

The aim of this section is to prove a criterion of intrinsic
formality.

Let $H$ be a graded $\CO$-algebra and let $\fg$ be the dg Lie algebra of
coderivations of the corresponding dg $\CO^{\perp}$-coalgebra
$(\mathbb{F}^*_{\CO^{\perp}}(H),Q)$. Since $\mathbb{F}^*_{\CO^{\perp}}(H)$
is cofree, any coderivation is uniquely defined by its composition
with the projection onto $H$. Therefore, $\fg$ considered
as a graded vector space, is isomorphic to 
$\Hom(\mathbb{F}^*_{\CO^{\perp}}(H),H)$. We denote 
$$\fg_{\ge 1}=\Hom(\oplus_{i\ge 2}\mathbb{F}^{*i}_{\CO^{\perp}}(H),H).$$
This is a dg Lie subalgebra of $\fg$.
\subsubsection{}
\begin{thm}{crit-i-f}
Suppose that the map $H^1(\fg_{\ge 1})\to H^1(\fg)$ is zero. 
Then $H$ is intrinsically formal.
\end{thm} 

\subsection{Proof of \Thm{crit-i-f}}

The following standard lemma results from the fact that $\CO_{\infty}$
is {\em cofibrant}.

\subsubsection{}
\begin{lem}{HPT}
Let $X$ be a $\CO_{\infty}$-algebra. There exists a $\CO_{\infty}$-algebra
structure on $H(X)$ so that $X$ and $H(X)$ are equivalent 
$\CO_{\infty}$-algebras (i.e., there exists a pair of quasi-isomorphisms
of  $\CO_{\infty}$-algebras $X\leftarrow F\to H(X)$).
\end{lem} 

\subsubsection{}
Let $H$ be a graded $\CO$-algebra. Let $X$ be a $\CO_{\infty}$-algebra
so that $H=H(X)$ as $\CO$-algebras. Choose a $\CO_{\infty}$-algebra
structure on $H$ guaranteed by \Lem{HPT}. One has $\CO_{\infty}(2)=\CO(2)$
and the $\CO$-algebra structure on $H$ is the restriction of the 
$\CO_{\infty}$-algebra structure. To fix a notation, let the collection
of maps
\begin{equation}
\label{o-infty-str}
Q_n:\mathbb{F}^{*n}_{\CO^{\perp}}(H)\to H[1],
\end{equation}
$n\ge 2$ define the said $\CO_{\infty}$-algebra structure on $H$.
The $\CO$-algebra structure on $H$ is given by the collection 
$\{Q^0_n\}$ with $Q^0_2=Q_2;\quad Q^0_i=0$ for $i>2$.

\subsubsection{}
\begin{lem}{}
Let $\lambda\in k$. Put $Q^{\lambda}_n=\lambda^{n-2}Q_n$.
The collection $\{Q^{\lambda}_n\}_{n\ge 1}$ defines a collection of 
$\CO_{\infty}$-algebra structures on $H$
parametrized by $\lambda\in k$. This gives the structure $\{Q_n\}$
for $\lambda=1$ and  $\{Q^0_n\}$ for $\lambda=0$.
\end{lem}
\begin{proof}

The only property we have to check to make sure that the collection 
$\{Q^{\lambda}_n\}$ defines a $\CO_{\infty}$-algebra structure, is
the identity looking like
$$d(Q^{\lambda}_n)=P_n(Q^{\lambda}_2,\ldots,Q^{\lambda}_{n-1})$$
where $P_n$ is a quadratic (non-commutative) polynomial.
Since $H$ has zero differential (this means $Q_1=0$ in our notation)
the left-hand side vanishes. The right hand side vanishes for $\lambda=1$
since the collection of $Q_i$ does define a $\CO_{\infty}$-action.
Since the polynomials $P_n$ are homogeneous, one has
$$P_n(Q^{\lambda}_2,\ldots,Q^{\lambda}_{n-1})=
\lambda^{n-1}P_n(Q_2,\ldots,Q_{n-1}).$$
This proves the claim.
\end{proof}

\subsubsection{}Put $C=(\mathbb{F}^*_{\CO^{\perp}}(H),Q^0)$.
This is a differential graded $\CO^{\perp}$-coalgebra. The
collection $\{Q^{\lambda}_n\}$ defines a $k[\lambda]$-linear differential
$Q^{\lambda}$ on the $\CO^{\perp}$-coalgebra $C[\lambda]$.
We wish to construct an isomorphism
$$ \theta:(C[\lambda],Q^0)\to (C[\lambda],Q^{\lambda})$$
which is identity modulo $\lambda$.

The isomorphism $\theta$  is uniquely defined by a collection of maps
$$\theta_n: \mathbb{F}^{*n}_{\CO^{\perp}}(H)\to H[\lambda]$$
with $\theta_1=\id_H$. We will be looking for $\theta$ satisfying the
following property.
\begin{equation}
\label{homog}
\theta_n=\phi_n\cdot\lambda^{n-1}\text{ for some }
\phi_n:\mathbb{F}^{*n}_{\CO^{\perp}}(H)\to H.
\end{equation}

An automorphism $\theta$ satisfying~(\ref{homog}) is constructed
in~\ref{constr-theta} below.
Then, tensoring $\theta$ by 
$k[\lambda]/(\lambda-1)$, we get an isomorphism of dg $\CO^{\perp}$-coalgebras 
$$\ol{\theta}:(C,Q^0)\iso (C,Q).$$
This will prove \Thm{crit-i-f}.

\subsubsection{}
Define an action of the multiplicative group $k^*$ on
$C[\lambda]$ by the formulas
$$ \mu * x=\mu^n\cdot x\text{ for }x\in\mathbb{F}^{*n}_{\CO^{\perp}}(H);
\quad \mu *\lambda=\mu\cdot\lambda.$$
The differentials $Q$ and $Q^{\lambda}$ have both degree $-1$ with respect to
this action:
$$  \mu * Q(\mu^{-1}* x)=\mu^{-1}\cdot Q(x);\quad
\mu * Q^{\lambda}(\mu^{-1}* x)=\mu^{-1}\cdot Q^{\lambda}(x).$$

The condition~(\ref{homog}) means that $\theta$ has degree zero with 
respect to the defined action of $k^*$.

\subsubsection{}
\label{constr-theta}
The map $\theta$ will be constructed by induction.

Suppose we have constructed an isomorphism
$$\theta: (C[\lambda]/(\lambda^n),Q^0)\iso (C[\lambda]/(\lambda^n),
Q^{\lambda})$$
satisfying the property $\theta_k=\phi_k\cdot\lambda^{k-1}$
for some $\phi_k:\mathbb{F}^{*k}_{\CO^{\perp}}(H)\to H$ for all $k$. 
This means in particular that $\theta_k=0$ for $k>n$. 

Our aim is to lift $\theta$ to a map 
$$\widetilde{\theta}: (C[\lambda]/(\lambda^{n+1}),Q)\iso 
(C[\lambda]/(\lambda^{n+1}),Q^{\lambda})$$
such that its components $\widetilde{\theta}_k$ satisfy the same property. 

First of all, we lift $\theta$ to the isomorphism
$$\theta':(C[\lambda]/(\lambda^{n+1}),Q')\iso 
(C[\lambda]/(\lambda^{n+1}),Q^{\lambda})$$
taking $\theta'_k=\theta_k$ for all $k$, where $Q'$ is some
differential uniquely defined by the above formula.  
The differential $Q'$ has also degree $-1$. Since $Q'$ coincides with $Q_0$
modulo $\lambda^n$, one has actually an equality $Q'_k=Q^0_k$
for $k\le n+1$ and $Q'_{n+2}=\lambda^n\cdot z$ for some
$z:\mathbb{F}^{*n+2}_{\CO^{\perp}}(H)\to H$.
One easily observes that the element $z$ considered as a derivation,
is a cycle. Therefore, there is a derivation $u\in\fg^0$, such that
$z=du$. This gives an isomorphism
$$\eta=\exp(\lambda^n\cdot u): 
(C[\lambda]/(\lambda^{n+1}),Q)\iso(C[\lambda]/(\lambda^{n+1}),Q')$$
which is identity modulo $\lambda^n$. The inductive step will
be acomplished if we are able to find an isomorphism
between $(C[\lambda]/(\lambda^{n+1}),Q)$ and $(C[\lambda]/(\lambda^{n+1}),Q')$ 
having degree zero.

The components $\eta_k$ of $\eta$ are divisible by $\lambda^n$ for $k>1$.
An easy calculation shows that the collection 
$\kappa_k:\mathbb{F}^{*k}_{\CO^{\perp}}(H)\to H$ 
given by the formulas
$$\kappa_1=\id_H,\ \kappa_{n+1}=\eta_{n+1},\ \kappa_i=0\text{ for } i\ne 1, 
n+1,$$
defines an isomorphism 
$$\kappa:(C[\lambda]/(\lambda^{n+1}),Q)\iso(C[\lambda]/(\lambda^{n+1}),Q').$$
The composition of $\kappa$ with $\theta'$ is the isomorphism 
$\widetilde{\theta}$ we were looking for.

The construction of isomorphism $\theta$ satisfying~(\ref{homog}), and, 
therefore, the proof of \Thm{crit-i-f}, is acomplished.

\section{Hochschild complex}
\label{hoch}

\subsection{Hochschild complex}
Let $A$ be an associative $k$-algebra. Its Hochschild complex $\CC:=C^*(A;A)$
has components defined by the formula
$$\CC^n=C^n(A;A)=\Hom(A^{\otimes n}, A),\quad n=0,1,\ldots$$
The graded vector space $\CC$ admits a $\Lie\{1\}$-algebra structure
which comes from the identification of $\CC[1]$ with the collection of
coderivations of the cofree coalgebra (with counit) cogenerated by $A[1]$.

An explicit formula for the Lie bracket is given in~\ref{expl-lie} below.

The multiplication
$\mu:A^{\otimes 2}\to A$ belongs to $\CC^2$; therefore the operator
$\ad \mu$ has degree $1$. An easy calculation shows that $(\ad \mu)^2=0$; 
$\CC$ endowed with the differential $\ad \mu$ becomes a dg $\Lie\{1\}$-algebra.

\subsection{$H(\CC)$ is a $\CG$-algebra} 
In order to prove \Thm{KFT} it would be enough to check that $H=H(C^*(A;A))$
is intrinsically formal as a Lie algebra. This, however, is not
true. Tamarkin's idea is to prove that $H$ becomes intrinsically formal when
it is considered as an algebra over an operad $\CG$ described below.

Define $m:\CC\otimes\CC\to\CC$ by the formula
\begin{equation}
\label{mult.}
 m(x\otimes y)=\mu\circ(x\boxtimes y)
\end{equation}
where $x\boxtimes y:A^{\otimes m+n}\to A^{\otimes 2}$ is defined to be
the tensor product of the maps $x:A^{\otimes m}\to A$ and 
$y:A^{\otimes n}\to A$.

The following lemma is due to M.~Gerstenhaber (1964).
\subsubsection{}
\begin{lem}{homol-is-g}The map $m$ induces a commutative associative 
multiplication on $H(\CC)$. The bracket on $H(\CC)$ is a derivation with 
respect to $m$.
\end{lem}

\subsubsection{}
\begin{defn}{gerst}Operad $\CG$ is the operad 
generated by the operations $m\in\CG(2)^0,\quad \ell\in\CG(2)^{-1}$ 
satisfying the following identities:

\begin{itemize}
\item $m$ is commutative associative
\item $\ell$ is Lie
\item $\ell$ is a derivation with respect to $m$.
\end{itemize}
\end{defn}

\Lem{homol-is-g} above means that the cohomology $H(\CC(A;A))$ admits
a natural $\CG$-algebra structure.

\subsubsection{}
The following construction assigns a $\CG$-algebra to any Lie algebra $\fg$.
Put $X=\mathbb{F}_{\Com}(\fg[-1])=\oplus_{i>0}S^i(\fg[-1])$. There is a unique 
$\Lie\{1\}$-algebra structure on $X$ extending that on $\fg[-1]$ such that 
$X$ becomes a $\CG$-algebra.
This is {\em a free $\CG$-algebra generated by a Lie algebra $\fg$}.

\subsubsection{}
\label{polyvectorfields}
There is a twisted (=sheaf) version of the above construction. Let $\fg$ be a 
Lie algebroid over a commutative algebra $A$. This mean that $\fg$ is a Lie 
algebra, $A$-module, and a map of Lie algebras and $A$-modules 
$\pi:\fg\to\Der(A,A)$ is given so that
$$[f,ag]=a[f,g]+\pi(f)(a)g$$
for $a\in A,\quad f,g\in\fg$.

Then a $\CG$-algebra structure on the $A$-symmetric algebra without unit
$S^{\ge 1}_A(\fg[-1])$ is naturally defined. If one defines $A=S^0_A(\fg[-1])$
to commute with $S^{\ge 1}_A(\fg[-1])$, one obtains a $\CG$-algebra
structure on the $A$-symmetric algebra $S_A(\fg[-1])$.

\subsection{Koszulity}

The operad $\CG$ is obviously quadratic. The quadratic dual cooperad
$\CG^{\perp}$ has as cogenerators elements 
$\widetilde{m},\quad\widetilde{\ell}$
of degrees $-1$ and $-2$ respectively. A simple calculation gives
\subsubsection{}
\begin{lem}{g-dual}
$$ \CG^{\perp}=\CG^*[2].$$
\end{lem}

One has the following important
\subsubsection{}
\begin{prop}{g-is-koszul}(\cite{gj})
$\CG$ is Koszul.
\end{prop}

For an easy proof of this fact see~\ref{pf-g-is-koszul}.

Recall that koszulity of $\CG$  means that the natural map (\ref{strictishom})
$$ \CGI\to\CG$$
is a quasi-isomorphism of operads. The operad $\CGI$
is the operad for {\em homotopy Gerstenhaber algebras}.

Deformation theory approach to the Formality Theorem is based on the following
version of {\em Deligne's conjecture}.

\subsubsection{}
\begin{thm}{dc}There is a 
structure of $\CGI$-algebra natural in $A$
on $C^*(A;A)$ inducing the described above $\CG$-algebra structure on 
$H(C^*(A;A))$.
\end{thm}

\Thm{dc} will be proven in Sections~\ref{bet} and \ref{two}. In this section we
will deduce Formality \Thm{KFT} from \Thm{dc}. 

\subsection{Calculation}

From now on $A$ is a smooth commutative $k$-algebra.
Our aim is to calculate the cohomology of $H:=H(C^*(A;A))$
and to make sure it vanishes when $A$ is a polynomial algebra.
This, together with \Thm{dc}, gives Formality Theorem. 

The following classical result of Hochschild-Kostant-Rosenberg describes
the cohomology of $C^*(A;A)$.
 
\subsubsection{}
\begin{lem}{HKR} $H=S_A(T_A[-1])$  where $T_A=\Der(A,A)$. 
The $\CG$-algebra structure on $H$ is defined as in~\ref{polyvectorfields}.
\end{lem}

Following \ref{crit-i-f}, we have to calculate the dg Lie algebra of 
coderivations  of the dg $\CG^{\perp}$-coalgebra 
$(\mathbb{F}^*_{\CG^{\perp}}(H),Q)$ corresponding to $H$.

\subsubsection{}
Note the following formula
\begin{equation}
\label{cofreegp=}
 \mathbb{F}^*_{\CG^{\perp}}(X)=\mathbb{F}^*_{\Com}
\left(\mathbb{F}^*_{\Lie}(X[1])[1]\right)[-2]
\end{equation}
which can be obtained using \ref{g-dual} from the formula
dual to the following
\begin{equation}
\label{freeg=}
\mathbb{F}_{\CG}(X)=\mathbb{F}_{\Com}(\mathbb{F}_{\Lie\{1\}}(X)).
\end{equation}

\subsubsection{}
\label{calc-1}
According to ~\ref{crit-i-f}, we have to calculate the 
map $H^1(\fg_{\ge 1})\to H^1(\fg)$ where 
$$\fg=\Coder(\mathbb{F}^*_{\CG^{\perp}}(H))=
\Hom(\mathbb{F}^*_{\CG^{\perp}}(H),H)=
\Hom(\mathbb{F}^*_{\Com}\left(\mathbb{F}^*_{\Lie}(H[1])[1]\right),H[2])$$
with the differential induced by the differential $Q$ of
$\mathbb{F}^*_{\CG^{\perp}}(H)$. 

The differential $Q$ of $\mathbb{F}^*_{\CG^{\perp}}(H)$ comes from the
map $\CG(2)\otimes H^{\otimes 2}\to H$ describing the $\CG$-algebra 
structure on $H$. Therefore, $Q=Q_m+Q_{\ell}$ where $Q_m$ is induced by the 
commutative multiplication $m:H\otimes H\to H$, and $Q_{\ell}$ is 
induced by the bracket $\ell: H\otimes H\to H[-1]$. Since the 
defining relations on operations $m$ and $\ell$ in $\CG$ are 
homogeneous, one necessarily has
$$ Q_m^2=Q_{\ell}^2=Q_mQ_{\ell}+Q_{\ell}Q_m=0.$$

The total diferential $Q$ on $\fg$ is also a sum of two differentials which
will be denoted by $Q_m$ and $Q_{\ell}$. 

Any cofree coalgebra is naturally graded --- see~(\ref{cofree}).
Formula~(\ref{cofreegp=}) gives rise to a bigrading on the cofree 
$\CG^{\perp}$-coalgebra $\mathbb{F}^*_{\CG^{\perp}}( H)$
in which the $(p,q)$-component consists of the elements of
$\Com$-degree $-p$ and total $\Lie$-degree $-q$.

This defines a bigrading on $\fg$
so that
$$ \fg^{p\bullet}=
\Hom(\mathbb{F}^{*1+p}_{\Com}(\mathbb{F}^{*\bullet}_{\Lie}( H[1])[1]),
 H[2]).$$

Note that
\begin{equation}
\fg_{\ge 1}=\bigoplus_{(p,q)\ne(0,0)}\fg^{pq}.
\end{equation}

The differentials $Q_m$ and $Q_{\ell}$ have degrees $(0,1)$ and $(1,0)$ with
respect to this bigrading and $\fg$ lives in the first quadrant. 
Therefore, one can use the spectral sequence argument to calculate 
the cohomology of $\fg$.

Let us calculate the first term $E_1^{pq}=H^{pq}(\fg,Q_m)$.
To keep track of the differential $Q_m$ in $\fg$ it is convenient to present
$$\fg^{0q}=\Hom(\mathbb{F}^{*1+q}_{\Lie}( H[1])[1], H[2])=
\Hom_{H}(\mathbb{F}^{*1+q}_{\Lie}( H[1])\otimes H[1], H[2])$$
and to identify $\mathbb{F}^{*}_{\Lie}( H[1])\otimes H$ with the
 homological Harrison complex $Z:=\Har_*(H,H)$.

Then one can see that $(\fg^{0\bullet},Q_m)$ coincides with 
$\Hom_{H}(Z[1], H[2])$
as a complex; moreover, for each $p$ one has
$$(\fg^{p\bullet},Q_m)=\Hom_{H}(S^{1+p}_{H}(Z[1]), H[2]).$$

\subsubsection{}
\label{calc-2}
The considerations above hold for every graded $\CG$-algebra $H$ with unit. 
Now we will use the fact that $ H=H(C^*(A;A))$ where $A$ is a smooth 
$k$-algebra.

Namely, according to~\Lem{HKR}, $H$ is smooth as a graded commutative 
algebra. Therefore, there is a natural isomorphism
$$ Z\iso\Omega[1]$$
where $\Omega=\Omega_{H/k}$ is the module of K\"ahler differentials.
This implies that
\begin{equation}
\label{e1} 
E_1^{pq}=\begin{cases}
                    \Hom_{H}(S^{1+p}_{H}(\Omega[2]), H[2]),
                              & q=0\\
                    0, &q\ne 0.
                \end{cases} 
\end{equation}

Let us calculate $\Omega_{H/k}$.
The sequence of smooth morphisms of graded commutative algebras
$$ k\to A\to H=S_A(T_A[-1])$$
gives rise to an isomorphism
\begin{equation}
  \label{Omega}
\Omega\iso H\otimes_A\Omega_{A/k}\oplus\Omega_{H/A}=
H\otimes_A(T_A[-1]\oplus T_A^*)=H\otimes_A\omega  
\end{equation}
where $\omega=T_A[-1]\oplus T_A^*$. 
Note that $\omega$ is a finitely generated graded projective $A$-module.

The only non-vanishing cohomology in~(\ref{e1}) can be rewritten as
\begin{equation}
\label{e1fin}
E_1^{p0}=S^{1+p}_A(\omega[-1])\otimes_A H[2]=S^{1+p}_H(\Omega[-1])[2]
\end{equation}
since $\omega[2]^*=\omega[-1]$.

Note that $E_1^{p0}$ embeds into 
$$\fg^{p0}=\Hom(\mathbb{F}_{\Com}^{*1+p}(H[2]),H[2])$$
and the differential $Q_{\ell}$ on the latter is defined by the Lie algebra
structure on $H[1]$.
This allows one to identify the the differential $Q_2$ on~(\ref{e1fin}) 
with the differential on the (shifted and truncated) de Rham complex of $H$.

\subsubsection{}
Suppose now that $A$ is a polynomial algebra over $k$. In this case
de Rham complex of $H$ is acyclic. Then the calculation in the previous
subsection gives a quasi-isomorphism $\fg\iso H/k[1]$. This implies that
$\fg$ has no cohomology coming from the cohomology of $\fg_{\ge 1}$.

\subsubsection{}
\begin{rem}{pf-g-is-koszul}
A calculation similar to the above proves that $\CG$ is Koszul.

In fact, according to~\ref{koszul-crit}, one has to check that for each
graded vector space $V$ the natural map
$$V\to(\mathbb{F}^*_{\CG^{\perp}}(\mathbb{F}_{\CG}(V),Q)$$
is a quasi-isomorphism.

Taking into account the formulas~(\ref{cofreegp=}) and~(\ref{freeg=})
and using, as in~\ref{calc-1}, the presentation of 
$\mathbb{F}^*_{\CG^{\perp}}(H)$ by a bicomplex, one easily
obtains the result. 
\end{rem} 

\subsubsection{}
Now \Thm{KFT} have been proven modulo \Thm{dc}. In the end of this section we
will describe the operad $\BI$ naturally acting on the Hochschild
complex $C(A;A)$ of any associative algebra $A$. This is the first step
in the proof of \Thm{dc} which is presented in Sections~\ref{bet} and
\ref{two}.

\subsection{Hochschild complex is a $\BI$-algebra}
We shall now describe an operad which acts naturally on the Hochschild
complex of any associative algebra. This operad is denoted $\BI$.
It has been invented by H.-I. Baues; its action on the Hochschild complex
was defined in~\cite{gj}.

\subsubsection{Notation}
In this subsection $A$ is any associative $k$-algebra and $\CC=C^*(A;A)$.
It is convenient to denote elements $f\in\CC^n$ as boxes having $n$ hands
and one leg like this:

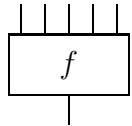
\begin{center}
\begin{picture}(10,3)
   \put(4,1){\makebox(2,1){$f$}}
   \put(4,1){\line(1,0){2}}
   \put(4,2){\line(1,0){2}}
   \put(4,1){\line(0,1){1}}
   \put(6,1){\line(0,1){1}}
   \put(5,1){\line(0,-1){0.5}}
   \put(4.2,2){\line(0,1){0.5}}
   \put(4.6,2){\line(0,1){0.5}}
   \put(5,2){\line(0,1){0.5}}
   \put(5.4,2){\line(0,1){0.5}}
   \put(5.8,2){\line(0,1){0.5}}
\end{picture}
\end{center}

\subsubsection{Basic operation} Let $f,g_1,\ldots,g_n\in\CC$.
Denote {\em the brace} $f\{g_1,\ldots,g_n\}$ by the following formula
\begin{equation}
\label{brace}
f\{g_1,\ldots,g_n\}=\sum_{\text{all possible insertions}}
\begin{picture}(10,4)
   \put(0,0){\makebox(10,1){$f$}}
   \put(0,0){\line(1,0){10}}
   \put(0,1){\line(1,0){10}}
   \put(0,0){\line(0,1){1}}
   \put(10,0){\line(0,1){1}}
   \put(5,0){\line(0,-1){0.5}}
   \put(0.3,1){\line(0,1){2}}
   \put(2.1,1){\line(0,1){2}}
   \put(2.4,1){\line(0,1){2}}
   \put(9.7,1){\line(0,1){2}}
   \put(9.4,1){\line(0,1){2}}
   \put(0.6,1.5){\makebox(1.2,1){$g_1$}}
   \put(2.7,1.5){\makebox(1.2,1){$g_2$}}
   \put(7.9,1.5){\makebox(1.2,1){$g_n$}}
   \put(0.6,1.5){\line(1,0){1.2}}
   \put(0.6,2.5){\line(1,0){1.2}}
   \put(2.7,1.5){\line(1,0){1.2}}
   \put(2.7,2.5){\line(1,0){1.2}}
   \put(7.9,1.5){\line(1,0){1.2}}
   \put(7.9,2.5){\line(1,0){1.2}}
   \put(0.6,1.5){\line(0,1){1}}
   \put(1.8,1.5){\line(0,1){1}}
   \put(2.7,1.5){\line(0,1){1}}
   \put(3.9,1.5){\line(0,1){1}}
   \put(7.9,1.5){\line(0,1){1}}
   \put(9.1,1.5){\line(0,1){1}}
   \put(0.9,2.5){\line(0,1){0.5}}
   \put(1.2,2.5){\line(0,1){0.5}}
   \put(1.5,2.5){\line(0,1){0.5}}
   \put(3,2.5){\line(0,1){0.5}}
   \put(3.3,2.5){\line(0,1){0.5}}
   \put(3.6,2.5){\line(0,1){0.5}}
   \put(8.2,2.5){\line(0,1){0.5}}
   \put(8.5,2.5){\line(0,1){0.5}}
   \put(8.8,2.5){\line(0,1){0.5}}
   \put(1.2,1){\line(0,1){0.5}}
   \put(3.3,1){\line(0,1){0.5}}
   \put(8.5,1){\line(0,1){0.5}}
\end{picture}
\end{equation}

Here the sum is taken over all possible order preserving insertions
of legs of $g_i$ into hands of $f$.

\subsubsection{}
\label{expl-lie}
The Lie bracket on $\CC[1]$ is given explicitly, in terms of  braces,
by the formula
$$ [f,g]=f\{g\}-(-1)^{|f||g|}g\{f\}.$$

\subsubsection{}
\begin{defn}{binfty} A $\BI$-algebra structure on a graded vector space $X$
is given by a structure of dg bialgebra on $\mathbb{F}_{\Ass}^*(X[1])$ so that
the coalgebra structure is the standard (cofree) one.
\end{defn}

Let us check that $\BI$-algebra structure is given by an operad
(as usual, it will be denoted by $\BI$).

The dg bialgebra structure on $\mathbb{F}_{\Ass}^*(X[1])$ is given by the 
following data.

\begin{itemize}

\item a differential $X[1]^{\otimes n}\to X[1]^{\otimes m}$ of degree $1$.
The differential is uniquely defined by its $m=1$ part. We denote its
$(n,1)$-components by $m_n:X[1]^{\otimes n}\to X[2]$ (or, what is the same,
$m_n:X^{\otimes n}\to X[2-n]$).

\item a multiplication $X[1]^{\otimes p}\otimes X[1]^{\otimes q}\to
X[1]^{\otimes r}$ of degree $0$ --- it is also uniquely defined by its 
$r=1$ part. We denote the collection of $r=1$ multiplications by
$$ m_{pq}:X[1]^{\otimes p}\otimes X[1]^{\otimes q}\to X[1]$$
or, what is the same,
$$ m_{pq}:X^{\otimes p}\otimes X^{\otimes q}\to X[1-p-q].$$
\end{itemize}

Therefore, the $\BI$-algebra structure is given by a collection of
operations $m_n,\ m_{pq}$ subject to some relations. This defines an
operad $\BI$ as the one generated by $m_n\in\BI(n)^{2-n}$ and 
$m_{pq}\in\BI(p+q)^{1-p-q}$ subject to some relations.

\subsubsection{}
{\bf WARNING.} The operad $\BI$ is not obtained in any sense from
a(ny) Koszul operad $\CB$. Getzler and Jones are responsible for this 
notation.

\subsubsection{Action of $\BI$ on $C^*(A;A)$}
\label{c-is-bi}
We have to define the action of the operations $m_n,\ m_{pq}$ on 
$\CC=C^*(A;A)$ and to check the compatibilities. Here it is.

\begin{itemize}
\item $m_1$ is the differential in $\CC$
\item $m_2$ is the multiplication $\mu$ defined by~(\ref{mult.})
\item $m_i=0$ for $i>2$
\item $m_{1k}(f\otimes g_1\otimes\ldots\otimes g_k)=f\{g_1,\ldots,g_k\}$
where the brace operations are defined by  formula~(\ref{brace})
\item $m_{kl}=0\text{ for } k>1.$
\end{itemize}
One can directly check that the collection of operations $m_n,\ m_{pq}$
defined above gives rise to a $\BI$-algebra structure on $\CC$.

\section{Between $\CG$ and $\CGI$}
\label{bet}

In this section we present an operad $\tB$ lying between 
$\CG$ and $\CGI$ --- it admits a pair of maps
$$\CGI\to\tB,\quad \tB\to\CG$$
so that  the composition is the obvious map $\CGI\to\CG$.

In the next section we will prove, using Etingof-Kazhdan theorem
on quantization of Lie bialgebras~\cite{ek}, that the operad $\tB$
is isomorphic to the operad $\BI$ acting on the Hochschild complex
of any associative algebra by~\ref{c-is-bi}. This will yield~\Thm{dc}
and, therefore, \Thm{KFT}.

\subsection{$\tB$-algebras}
A $\tB$-algebra structure on a graded vector space $X$ is a dg Lie
bialgebra structure on $\mathbb{F}^*_{\Lie}(X[1])$ extending the standard free 
Lie coalgebra structure.

The Lie bracket on a Lie bialgebra  $\mathbb{F}^*_{\Lie}(X[1])$ is defined 
by its corestriction to the cogenerators $X[1]$. Therefore, it is given by
a collection of maps
$$ \ell_{mn}:\mathbb{F}^{*m}_{\Lie}(X[1])\otimes \mathbb{F}^{*n}_{\Lie}(X[1])
\to X[1]$$
satisfying a collection of quadratic identities.
The differential on  $\mathbb{F}^*_{\Lie}(X[1])$ is also defined by its 
corestriction to the cogenerators. This amounts to a collection
$$ d_n: \mathbb{F}^{*n}_{\Lie}(X[1])\to X[2]$$
satisfying some more quadratic identities --- the one saying that $d^2=0$
and the other that the $d$ is the derivation of the Lie algebra structure 
given by $\ell_{mn}$.

In particular, one has $d_1^2=0$ and this endows $X$ with a structure of 
complex. The  obvious maps  $X[1]\to \mathbb{F}^*_{\Lie}(X[1])\to X[1]$ are 
maps of complexes.

Since the $\tB$-structure on $X$ is given by a collection of operations
subject to some relations, there is an operad in the category of 
complexes which will be called in the sequel $\tB$ such that 
$\tB$-algebras are just algebras over $\tB$. 

\subsection{}
A map $\CO\to\CO'$ of operad endows a $\CO'$-algebra with a
canonical $\CO$-algebra structure. The converse is also obviously true --- 
in order to define a map of operads it is enough to endow any $\CO'$-algebra
with a canonical $\CO$-algebra structure.

Let us construct a map $\tB\to\CG$. For this we have to define 
canonically a $\tB$-algebra structure on each $\CG$-algebra $X$.
Recall that a $\CG$-algebra $X$ is endowed with a commutative
multiplication $m:X^{\otimes 2}\to X$ and a Lie bracket
$l:X[1]^{\otimes 2}\to X[1]$. The Harrison complex of the commutative
algebra $(X,m)$ is given by a differential on $\mathbb{F}^*_{\Lie}(X[1])$.
The Lie algebra structure on $X[1]$ can be uniquely extended to  
 $\mathbb{F}^*_{\Lie}(X[1])$ to get a Lie bialgebra. The Harrison differential
will be a derivation with respect to the Lie algebra structure, so this
construction defines a dg Lie bialgebra structure on $\mathbb{F}^*_{\Lie}(X[1])$.

The construction is obviously canonical and yields a morphism of operads
$\tB\to\CG$.

\subsection{} Let now construct a map $\CGI\to\tB$.

Let $X$ be a $\tB$-algebra. This means that a dg Lie bialgebra structure
on $\fg=\mathbb{F}^*_{\Lie}(X[1])$ is given. In particular, $\fg$ is a dg 
Lie algebra and this defines a differential on $\mathbb{F}^*_{\Com}(\fg[1])$.
 The latter complex is by formula~(\ref{cofreegp=}) just 
$\mathbb{F}^*_{\CG^{\perp}}(X)[2]$. 
Differential on it gives a $\CGI$-structure on $X$.

Thus the map $\CGI\to\tB$ is constructed.

\section{Equivalence of $\tB$ with $\BI$}
\label{two}

In this section we prove that the operads $\tB$ and $\BI$ are isomorphic.
Note that the isomorphism is obtained using Etingof-Kazhdan theorems~\cite{ek}
on quantization of Lie bialgebras. The isomorphism will depend on 
the choice of associator, as in~\cite{ek}.

\subsection{} For some technical reasons, it is more convenient to use
coalgebras over $\tB$ and $\BI$ instead of algebras. Our aim is to prove that
any $\tB$-coalgebra admits a natural $\BI$-coalgebra structure and vice versa.

Note

\subsubsection{}
\begin{lem}{tb-c}
$\tB$-coagebra structure on a graded vector space $X$ is given by a
structure of dg Lie bialgebra on
$$ \widehat{\mathbb{F}}_{\Lie}(X[1])=\prod_{n=0}^{\infty}
\mathbb{F}^n_{\Lie}(X[1]).$$
\end{lem}

\subsubsection{}
\begin{lem}{bi-c}
$\BI$-coagebra structure on a graded vector space $X$ is given by a
structure of dg bialgebra on
$$ \widehat{\mathbb{F}}_{\Ass}(X[1])=\prod_{n=0}^{\infty}X[1]^{\otimes n}.$$
\end{lem}

Now we wish to use~\cite{ek} in order to pass from one structure
above to the other. The idea is the following. One can interpret
completions of free algebras $\mathbb{F}_{\Lie}(V)$ and 
$\mathbb{F}_{\Ass}(V)$ as equivariant $k[[h]]$-algebras
 $\mathbb{F}_{\Lie}(V)[[h]]$ and $\mathbb{F}_{\Lie}(V)[[h]]$.
This is the situation Etingof-Kazhdan theory applies.

\subsection{Etingof-Kazhdan theory}
Let $\Locc(k)$ be the category of local complete $k$-algebras with residue 
field $k$. 

Let $\CA$ be an abelian $k$-linear tensor category. 
For each $R\in\Locc(k)$ we denote by $\CA(R)$
the category with the same objects as $\CA$ and with the morphisms defined by
the formula
$$ \Hom_{\CA(R)}(X,Y)=\Hom_{\CA}(X,Y)\otimes R.$$

The object of $\CA(R)$ corresponding to an object $X\in\CA$ is denoted
$X_R$. In the other direction, for $Y\in\CA(R)$ we write
$\ol{Y}$ for the corresponding object of $\CA$. The assignments 
$X\mapsto X_R$ and $Y\mapsto \ol{Y}$ define a pair of functors
between $\CA$ and $\CA(R)$.

Let $\LBA(R)$ be the category of Lie bialgebras $(\fg,[\quad],\delta)$ in 
$\CA(R)$ whose  cobracket $\delta$ vanishes modulo the maximal ideal $\fm$ 
of $R$. 
Let $\HA$ denote the category of Hopf algebras in $\CA_R$ whose reduction 
modulo $\fm$ is isomorphic to the enveloping algebra of a Lie algebra
in $\CA$.

The following theorem can be found in~\cite{ek}.

\subsubsection{}
\begin{thm}{q-deq}
There is an equivalence of categories
\begin{equation}
\label{quant}
Q:\LBA(R)\to\HA(R)
\end{equation}
satisfying the following properties (see also explanations below)
\begin{itemize}
\item[1.] $\ol{Q(\fg)}=U(\ol{\fg})$ 
\item[2.] $\delta:\fg\to\fg\otimes\fg$ measures the deviation of the 
coproduct in $Q(\fg)$ from being cocommutative.
\item[3.] $Q$ is given by universal formulas.
\end{itemize}
\end{thm}

\subsubsection{}
The second property of the functor $Q$ mentioned in~\Thm{q-deq} means the
following. Denote by $i:\fg\to Q(\fg)$ the image under the functor $\otimes R$
of the  obvious embedding $\ol{\fg}\to U(\ol{\fg})$. The property $2$ claims 
that the map from $\fg$ to $Q(\fg)\otimes Q(\fg)$ given by the difference
$$ (\Delta-\Delta')\circ i - (i\otimes i)\circ \delta$$
vanishes modulo $\fm^2$.

Here $\Delta$ is the coproduct in $Q(\fg)$ and $\Delta'$ is the coproduct 
composed with the commutativity constraint.

\subsubsection{}
The third property means the following. As an object of $\CA(R)$,
$Q(\fg)$ is just the symmetric algebra $S(\fg)=\oplus S^n(\fg)$.
Therefore, the Hopf algebra structure on $Q(\fg)$ is given by a collection
of maps $m_{pqr}: S^p(\fg)\otimes S^q(\fg)\to S^r(\fg)$ and 
$\Delta_{pqr}: S^p(\fg)\to S^q(\fg)\otimes S^r(\fg)$. Universality condition
means that the maps $m_{pqr},\quad\Delta_{pqr} $ are described as universal
polynomials on the bracket and cobracket in $\fg$.

\subsubsection{} It is convenient to define $\LBA$ to be the category of
pairs $(R,\fg)$ where $R\in\Locc(k)$ and $\fg\in\LBA(R)$. In the same
fashion one defines the category $\HA$. Since the functors 
$Q:\LBA(R)\to\HA(R)$ are given by the universal formulas, they form a functor
$Q:\LBA\to\HA$ which is also an equivalence of categories. Reduction 
modulo the maximal ideal defines a commutative diagram of functors
$$
\begin{CD}
\LBA@>>Q> \HA \\
@VVV @VVV \\
\Liealg(k)@>>U>\HA(k)
\end{CD}
$$
where $\Liealg(k)$ is the category of Lie algebras over $k$ and $U$ is the
enveloping algebra functor.

\subsubsection{}
\label{k*-actions}
The equivalence of categories $Q:\LBA\to\HA$ gives rise to
an equivalence $Q^G:\LBA^G\to\HA^G$ between the categories of objects
endowed with a $G$-action, $G$ being a group.

Let now $\CA$ be the category of complexes of $k$-modules.
Let $G=k^*$ be the multiplicative group, $R=k[[h]]$. Let $k^*$ act on $R$
by the formula $\lambda(h)=\lambda^{-1}h$. Let $V\in\CA$.
Let $k^*$ act on $V$ by the formula 
$\lambda(v)=\lambda\cdot v$. This action extends to a $k^*$-action
on $\mathbb{F}_{\Lie}(V)$ and $\mathbb{F}_{\Ass}(V)$, as well as to
an action on $\mathbb{F}_{\Lie}(V)[[h]]$ and $\mathbb{F}_{\Ass}(V)[[h]]$.

\Thm{q-deq} implies the following

\subsubsection{}
\begin{cor}{equivariant-qdq}
The functor $Q$ establishes an equivalence between the following categories:
\begin{itemize}
\item[1.] Lie bialgebras $(R,\fg)\in\LBA$ endowed with a $k^*$-action 
compatible with the specified above $k^*$-action on $R=k[[h]]$ 
and on $\ol{\fg}=\mathbb{F}_{\Lie}(V)$.

\item[2.] Associative bialgebras $(R,H)\in\HA$ endowed with a $k^*$-action 
compatible with the specified above action on $R=k[[h]]$ 
and on $\ol{H}=\mathbb{F}_{\Ass}(V)$.
\end{itemize}
\end{cor}

\subsection{}
\begin{thm}{b=b}
There exists an isomorphism between the operads $\tB$ and $\BI$.
\end{thm}  

\Thm{b=b} is proven in~\ref{bop}--\ref{eop} below.

\subsubsection{}
\label{bop}
Put $\fg=\mathbb{F}_{\Lie}(V)$. A $k[[h]]$-Lie bialgebra structure
on $\fg[[h]]$ is given by a collection of maps
$$ \delta_{pq}^r:V\to\mathbb{F}^p_{\Lie}(V)\otimes \mathbb{F}^q_{\Lie}(V)$$
such that the cobracket $\delta:\fg[[h]]\to\fg[[h]\otimes\fg[[h]]$
restricted to $\fg$ is given by the formula
$$\fg=\sum_{p,q,r}\delta_{pq}^r\cdot h^r.$$

Define an action of $k^*$ on $\fg[[h]]$ as in~\ref{k*-actions}.
The cobracket $\delta$ of $\fg[[h]]$ is equivariant if and only if
it satisfies the property
\begin{equation}
\label{eq-lie}
\delta_{pq}^r=0\text{ for }r\ne p+q-1.
\end{equation}

One can easily identify dg Lie bialgebra structures on $\fg[[h]]$ satisfying
(\ref{eq-lie}) with dg Lie bialgebra structures on $\widehat{\fg}$.

\subsubsection{}
Similarly, dg bialgebra structures on $\widehat{\mathbb{F}}_{\Ass}(V)$
can be identified with equivariant bialgebra structures on 
$\mathbb{F}_{\Ass}(V)[[h]]$.

\subsubsection{} 
\label{eop}
We use~\Cor{equivariant-qdq} of the equivalence $Q$ from 
Etingof-Kazhdan \Thm{q-deq}.

Let $X$ be a complex, $V=X[1]$. $\tB$-coalgebra
structure  on $X$ is given by a structure of dg Lie bialgebra on
$\widehat{\fg}=\widehat{\mathbb{F}}_{\Lie}(V)$
which is the same as an equivariant dg Lie bialgebra structure on $\fg[[h]]$.

According to~\Cor{equivariant-qdq}, this defines canonically an
equivariant dg Hopf algebra $(H,m,\Delta)\in\HA$.

The canonical map $i:V[[h]]\to H$ given by the composition
$$i:V\to\fg\to\ol{H}$$
induces an algebra homomorphism 
$F(i):\mathbb{F}_{\Ass}(V)[[h]]\to H$. It is
isomorphism since its reduction modulo $h$ $\ol{F(i)}$ is the identity map.
This defines canonically an equivariant bialgebra structure on 
$\mathbb{F}_{\Ass}(V)[[h]]$ which is the same as a dg bialgebra structure
on $\widehat{\mathbb{F}}_{\Ass}(V)$.
Theorem is proven.

\end{document}